\documentclass[11pt]{article}
\usepackage{amsmath}
\usepackage{amssymb}

\newtheorem{theorem}{Theorem}[section]
\newtheorem{lemma}{Lemma}[section]
\newtheorem{corollary}{Corollary}[section]

\begin{document}

\title{Asymptotic structure of viscous incompressible flow around a rotating body,
with nonvanishing flow field at infinity.
}

\author {
Paul Deuring\footnote{Univ Lille
Nord de France, 59000 Lille, France; ULCO, LMPA, 62228 Calais c\'edex, France.}
,\
Stanislav Kra\v cmar\footnote{
Department of Technical Mathematics, Czech Technical University,
Karlovo n\'{a}m. 13, 121 35 Prague 2, Czech Republic}
,\  \v S\'arka Ne\v casov\' a\footnote{Mathematical Institute
of the Academy of Sciences of the Czech Republic, \v Zitn\' a~25,
115 67 Praha 1, Czech Republic}}
\date{ }
\maketitle

\vspace{2ex}
\begin{abstract}
We consider weak (''Leray'') solutions to the stationary Navier-Stokes system with
Oseen and rotational terms, in an exterior domain. It is shown the velocity may be split
into a constant times the first column of the fundamental solution of the Oseen system,
plus a remainder term decaying pointwise near infinity at a rate which is higher than the
decay rate of the Oseen tensor. This result improves the theory by M. Kyed,
Asymptotic profile of a linearized flow past a rotating body,
Q. Appl. Math. 71 (2013), 489-500.



\end{abstract}

\section{Introduction}

Let $\mbox{$\mathfrak D$} \subset \mathbb{R}^3 $ be an open bounded set. Suppose this
set describes a rigid body moving with constant nonzero
translational and angular velocity in an incompressible viscous fluid.
Then the flow
aroung this body with respect to a frame attached to this body is governed by the following set
of non-dimensional equations (see \cite{G2}),
\begin{eqnarray} \label{10} &&
-\Delta u + \tau \, \partial _1 u  + \tau \, (u \cdot \nabla )u - (\omega \times x) \cdot \nabla  u +
\omega \times u + \nabla \pi = f,
\;\;
\mbox{div}\, u =0,
\end{eqnarray}
in the exterior domain
$
\overline{ \mbox{$\mathfrak D$} }^c:=
\mathbb{R}^3 \backslash \overline{ \mbox{$\mathfrak D$} },
$
supplemented by a decay condition at infinity,
\begin{eqnarray} \label{20}
u(x)\to 0\;\; \mbox{for}\;\; |x|\to \infty ,
\end{eqnarray}
and suitable boundary conditions on $\partial \mbox{$\mathfrak D$} $.
These latter conditions need not be specified here because they are not relevant in
the context of the work at hand. In (\ref{10}) and (\ref{20}), the functions
$u: \overline{ \mbox{$\mathfrak D$} }^c \mapsto
\mathbb{R}^3 $ and $\pi : \overline{ \mbox{$\mathfrak D$} }^c \mapsto
\mathbb{R}$
are the unknown relative velocity and pressure field of the fluid, respectively,
whereas the function
$f: \overline{ \mbox{$\mathfrak D$} }^c \mapsto
\mathbb{R}^3 $ stands for a prescribed volume force acting on the fluid.
The vector $\tau \, (-1,0,0)$
represent the uniform velocity of the flow at infinity or the velocity
of the body, depending
on the physical situation, and
$ \omega := \varrho \cdot (1,0,0)$ corresponds to the constant angular velocity of the body.
In particular the vectors of translational and angular velocity are parallel. From
a physical point of view this assumption is natural for a steady flow.
The parameters  $\tau \in (0, \infty ) $ and $ \varrho \in \mathbb{R} \backslash\{0\} $
are dimensionless quantities that can be identified with the Reynolds and Taylor number,
respectively. They will be considered as fixed, like the domain \mbox{$\mathfrak D$}.

We are interested in ``Leray solutions'' of (\ref{10}), (\ref{20}), that is, weak solutions characterized
by the conditions $u \in L^6( \overline{ \mbox{$\mathfrak D$} }^c)^3\cap
W^{1,1}_{loc}( \overline{ \mbox{$\mathfrak D$} }^c)^3,\; \nabla u  \in L^2( \overline{ \mbox{$\mathfrak D$} }^c)^9$
and $\pi \in L^2_{loc}( \overline{ \mbox{$\mathfrak D$} }^c)$.
The relation $u \in L^6( \overline{ \mbox{$\mathfrak D$} }^c)^3$ means that (\ref{20}) is verified in a weak
way.
Such solutions exist for data of arbitrary size
if some smoothness of the boundary of \mbox{$\mathfrak D$} is required, and if suitable regularity
conditions are imposed on $f$ and on the data on the boundary (\cite[Theorem IX.3.1]{Galdi-neu}).
It is well known by now (\cite{GK1}, \cite{DKN5}) that the velocity part $u$ of a Leray solution $(u,\pi )$
to (\ref{10}), (\ref{20}) decays for $|x|\to \infty $ as expressed by the estimates
\begin{eqnarray} \label{40}
|u(x)|\le C\, \bigl(\, |x|\, s(x) \,\bigr) ^{-1},
\quad
| \nabla u(x)|\le C\, \bigl(\, |x|\, s(x) \,\bigr) ^{-3/2}
\end{eqnarray}
for $x \in \mathbb{R}^3 $ with $|x|$ sufficiently large, where $s(x):=1+|x|-x_1\; (x \in \mathbb{R}^3 )$ and
$C>0$ a constant independent of $x$. The factor $s(x)$ may be considered as a mathematical manifestation
of the wake extending downstream behind a body moving in a viscous fluid.
In view of (\ref{40}), it is natural to ask how an asymptotic expansion of $u$ for $|x|\to \infty $
might look like. As far as we  know, up to now there are two answers to this question. The first is due
to Kyed \cite{K2}, who showed that
\begin{eqnarray} \label{60}
u_j(x)= \gamma \, E_{j1}(x) + R_j(x),
\quad
\partial _lu_j(x)= \gamma \, \partial _lE_{j1}(x) + S_{jl}(x)
\;\;
(x \in \overline{ \mbox{$\mathfrak D$} }^c,\; 1\le j,l\le 3),
\end{eqnarray}
where $E: \mathbb{R}^3 \backslash\{0\} \mapsto \mathbb{R} ^4 \times \mathbb{R}^3 $
denotes a fundamental solution to the Oseen system
\begin{eqnarray} \label{50}
- \Delta v+\tau \, \partial _1v+ \nabla \varrho = f,\;\;
\mbox{div}\, v =0
\quad \mbox{in}\;\; \mathbb{R}^3 .
\end{eqnarray}
The definition of the function $E$ is
stated in Section 2.
As becomes apparent from this definition, the term $E_{j1}(x)$ may be expressed
explicitly  in terms of elementary functions. The coefficient $\gamma $ is also given explicitly,
its definition involving the Cauchy stress tensor. The remainder terms $R$ and $S$ are characterized
by the relations $R \in L^q( \overline{ \mbox{$\mathfrak D$}  }^c)^3$ for $ q \in (4/3,\, \infty )$,
$S \in L^q( \overline{ \mbox{$\mathfrak D$}  }^c)^3$ for $ q \in (1, \infty )$.
Since it is known from \cite[Section VII.3]{Ga1} that
$E_{j1}|B_r^c \notin L^q(B_r^c)$ for $r>0,\; q \in [1,2]$, and
$\partial _lE_{j1}|B_r^c \notin L^q(B_r^c)$ for $r>0,\; q \in [1,\,4/3],\; j,l \in \{1,\, 2,\, 3\} $,
the function $R$ decays faster than $E_{j1}$, and $S_{jl}$ faster than $\partial _lE_{j1}$, in the
sense of $L^q$-integrability.
Thus the equations in (\ref{60}) may in fact be considered as asymptotic expansions of $u$ and
$\nabla u$, respectively. The theory in \cite{K2} is valid under the assumption that $u$ verifies
the boundary conditions
\begin{eqnarray} \label{61}
u(x)= e_1+ ( \omega \times x) \quad \mbox{for}\;\; x \in \partial \mbox{$\mathfrak D$} .
\end{eqnarray}
Reference \cite{K2} does not deal with pointwise decay of $R$ and $S$, nor does it indicate whether $S= \nabla R$.
The second answer to the question of how $u$ may be expanded asymptotically for $|x|\to \infty $ is
given in reference \cite{DKN7}, which states that for $x \in \overline{ B_{S_1}}^c,\;  1\le j\le 3,$
\begin{eqnarray} \label{70}
u_j(x)
=
\sum_{k = 1}^3 \beta _k\, Z_{jk}(x,0)
+ \Bigl( \int_{ \partial \Omega }u \cdot n\, do_x \Bigr) \, x_j\,(4\, \pi\, |x| ^3) ^{-1}
+ \mbox{$\mathfrak F$} _j(x) .
\end{eqnarray}
Here $S_1$ is a sufficiently large positive real number, $(Z _{jk} )_{1\le j,k\le 3}$ is the tensor
velocity part of the fundamental solution constructed by Guenther, Thomann \cite{GT} for the
linearization
\begin{eqnarray} \label{15}
-\Delta v + \tau \, \partial _1 v   - (\omega \times x) \cdot \nabla  v +
\omega \times v + \nabla \varrho  = f,
\;\;
\mbox{div}\, v =0
\end{eqnarray}
of (\ref{10}) (see Section 2 for the definition of $Z$), and \mbox{$\mathfrak F$} is a function
from $C^1( \overline{ B_{S_1}}^c)^3 $ given explicitly in terms of $Z,\, u$ and $\pi$ (Theorem \ref{theorem2.10}).
As is shown in \cite{DKN7}, this function \mbox{$\mathfrak F$} decays pointwise, in the sense that
\begin{eqnarray*}
\lim_{|x|\to \infty }| \partial ^{\alpha }\mbox{$\mathfrak F$} (x)| = O \bigl(\, (|x|\, s(x) )^{-3/2- |\alpha |/2}\,
\ln (2+|x|) \,\bigr)
\quad \mbox{for}\;\;
\alpha \in \mathbb{N} _0^3\; \mbox{with}\; | \alpha |\le 1.
\end{eqnarray*}
It is known from \cite[Theorem 2.19]{DKN2} -- and restated below in Corollary \ref{corollary2.120} -- that
\begin{eqnarray} \label{99}
\lim_{|x|\to \infty }| \partial ^{\alpha }Z (x,0)| = O \bigl(\, (|x|\, s(x) )^{-1- |\alpha |/2}\,\bigr)
\quad
(\alpha \in \mathbb{N} _0^3\; \mbox{with}\; | \alpha |\le 1).
\end{eqnarray}
So, if the decay rate in (\ref{99}) is sharp,
equation (\ref{70}) may be considered as an asymptotic expansion in the usual sense: the remainder exhibits
a faster pointwise decay than the leading term. The coefficients $\beta _1,\, \beta _2,\, \beta _3$
in (\ref{70}) are given explicitly in terms of $u,\, \pi$ and $f$. The theory in \cite{DKN7} does not
impose any boundary condition on $u$ or $\pi$. However, since the definition of the term $Z(x,0)$ involves an
integral over $(0, \infty )$, the leading term $\sum_{k = 1}^3 \beta _k\, Z_{jk}(x,0)$ in (\ref{70}) is not
as explicit as one would like it to be. More details on the theory from \cite{DKN7} may be found in
Theorem \ref{theorem2.10} below, where the main result from \cite{DKN7} is restated.

In the work at hand, we show that $Z_{j1}(x,0)= E_{j1}(x)$ for $x \in \mathbb{R}^3 \backslash\{0\} ,\;
1\le j\le 3$, and $\lim_{|x|\to \infty }| \partial ^{\alpha }_xZ _{jk} (x,0)|=O \bigl(\, (|x|\, s(x))^{-3/2-| \alpha |/2}
\,\bigr) $ for $1\le j\le 3,\; k \in \{2,\, 3\} $ (Corollary \ref{corollary3.50}, Theorem \ref{theorem4.10}).
Thus, setting
\begin{eqnarray} \label{ZU31}
\mbox{$\mathfrak G$} _j(x):=
\sum_{k = 2}^3 \beta _k\, Z_{jk}(x,0) + \mbox{$\mathfrak F$} _j(x)
\quad
(x \in \overline{ B_{S_1}}^c,\; 1\le j\le 3),
\end{eqnarray}
we may deduce from (\ref{70}) that
\begin{eqnarray} \label{80} \hspace{-1.5em}
u_j(x)
=
\beta _1\, E_{j1}(x)
+ \Bigl( \int_{ \partial \Omega }u \cdot n\, do_x \Bigr) \, x_j\, (4\, \pi\, |x| ^3) ^{-1}
+ \mbox{$\mathfrak G$} _j(x)
\;\;\;
(x \in \overline{ B_{S_1}}^c,\; 1\le j\le 3)
\end{eqnarray}
and
\begin{eqnarray} \label{90}
\lim_{|x|\to \infty }| \partial ^{\alpha }\mbox{$\mathfrak G$} (x)| = O \bigl(\, (|x|\, s(x) )^{-3/2- |\alpha |/2}\,
\ln (2+|x|) \,\bigr)
\quad \mbox{for}\;\;
\alpha \in \mathbb{N} _0^3\; \mbox{with}\; | \alpha |\le 1
\end{eqnarray}
(Theorem \ref{theorem2.20}, Corollary \ref{corollary2.10}).
If we compare how the coefficient $\gamma $ from (\ref{60}) is defined in \cite{K2},
and the coefficient $\beta _1$ from (\ref{80}) in \cite{DKN7} (see Theorem \ref{theorem2.10} below),
and if we take account of the boundary condition (\ref{61}) satisfied by $u$ in \cite{K2}, we see
that $\gamma $ and $\beta _1$ coincide. Thus the relations in (\ref{80}) and (\ref{90}) provide
a synthesis of the theories in \cite{K2} and \cite{DKN7}: the leading terms in (\ref{60}) and (\ref{80})
are identical, and the remainder in (\ref{80}) decays pointwise for $|x|\to \infty $, its rate of decay being
$|x| ^{-2-| \alpha |} +(|x|\, s(x) )^{-3/2- |\alpha |/2}\, \ln (2+|x|) $.
It is shown in \cite{KNP} -- and restated below in Theorem \ref{theorem2.70} -- that
\begin{eqnarray} \label{100}
\lim_{|x|\to \infty }\partial ^{\alpha }E_{j1}(x) = O\bigl(\, (|x|\, s(x) )^{-1- |\alpha |/2} \,\bigr)
\quad \mbox{for}\;\;
\alpha \in \mathbb{N} _0^3\; \mbox{with}\; | \alpha |\le 1
\end{eqnarray}
(\cite[(1.14)]{KNP}). The theory in \cite{KNP} additionally yields that the decay rate
$O(|x|\, s(x) )^{-1- |\alpha |/2}$ in (\ref{100}) is sharp. Therefore it follows from (\ref{90}) that
equation (\ref{80}) is in fact an asymptotic expansion of $u_j(x)$ for $|x|\to \infty $,
with the remainder vanishing faster for large values of $|x|$ than the leading term
$\beta _1\, E_{j1}(x) $. It further follows that the decay rates of $u$ and $\nabla u$
given by (\ref{40}) are sharp, too.
The reader may wish to check on the basis of the theory in \cite{KNP} whether some part of $E_{j1}(x)$
may be split off and put into the remainder term. (This is possible for $\nabla E_{j1}(x) $ but
not for $E_{j1}(x)$ if a decay rate as in (\ref{90}) is to be maintained.)

We further remark that in the case of a rigid body which only rotates but does not translate,
more detailed asymptotic expansions are available (\cite{FH092} -- \cite{FH2}).
Any reader interested in further results on the asymptotic behaviour of viscous incompressible flow
around rotating bodies is referred to
\cite{DKN1} -- \cite{DKN4}, \cite{DKN6}, \cite{Fa2} -- \cite{FH}, \cite{FHM} -- \cite{FN},
\cite{GaT}, \cite{Galdi-neu}, \cite{GK2} -- \cite{KNP-JAP}, \cite{KrPe1} -- \cite{K3}, \cite{Ne2} -- \cite{NS},
\cite{GT}.

\section{Notation. Definition of fundamental solutions. Auxiliary results.}
\setcounter{equation}{0}

By $| \;|$ we denote the Euclidean norm in $\mathbb{R} ^3$ and
the length $\alpha _1+ \alpha _2+ \alpha _3$ of a multiindex $\alpha \in \mathbb{N} _0^3$.
Put $e_1:=(1,0,0).$
For $r>0$, we set $B_r:=\{y \in \mathbb{R}^3 \: :\: |y|< r\}.$
If $A \subset \mathbb{R}^3 $, we put
$A^c:= \mathbb{R}^3 \backslash A$. Recall the abbreviation $s(x):= 1+|x|-x_1\; (x \in \mathbb{R}^3) $
introduced in Section 1.

If $A \subset \mathbb{R}^3 $ is open, $p \in [1, \infty )$ and $ k \in \mathbb{N}  $,
we write  $ W^{k,p}(A)$ for the usual Sobolev space of order $k$ and exponent $p$.
If $B \subset \mathbb{R}^3 $ is again an open set,
we define $ L^p_{loc}(B),\; W^{k,p}_{loc}(B)$ as the set of all functions $v: B \mapsto \mathbb{R} $
such that $v|_{U} \in L^{p}(U)$ and $v|_{U} \in W^{k,p}(U)$, respectively,
for any open bounded set $U \subset \mathbb{R}^3 $ with
$ \overline{ U} \subset B$.
We write $\mbox{$\mathfrak S$} ( \mathbb{R}^3 )$  for the usual space of rapidly decreasing functions in
$ \mathbb{R}^3 $; see \cite[p. 138]{Neri} for example.
For the Fourier transform $\widehat{g}$ of a function $g \in L^1( \mathbb{R}^3 )$, we choose the definition
$
\widehat{g}(\xi):=(2 \hspace{1pt}  \pi ) ^{-3/2} \hspace{1pt}
\int_{ \mathbb{R}^3 }e^{-i \hspace{1pt}  \xi \hspace{1pt}  x}
\hspace{1pt}  g(x)\; dx
\; (\xi \in \mathbb{R}^3 ).
$
This fixes the definition of the Fourier transform of a tempered distribution as well.

The numbers $\tau \in (0, \infty ) ,\; \varrho \in \mathbb{R} \backslash\{0\} $ introduced in Section 1
will be kept fixed throughout. Recall that the vector $\omega $ is given by $\omega:= r \cdot e_1.$
We introduce a matrix $ \Omega \in \mathbb{R}^{3 \times 3}  $ by setting
\begin{eqnarray*}
\Omega :=
\varrho \hspace{1pt}    \left( \begin{array}{rrr}
0 & 0  & 0  \\
0& 0 & - 1\\
0 & 1 & 0
\end{array}
\right).
\end{eqnarray*}
Note that $ \omega   \times x = \Omega \cdot x$ for $x \in \mathbb{R}^3 $.
We write \mbox{$\mathfrak C$} for positive constants that may depend on $\tau $ or $\varrho $.
Constants additionally depending on parameters $\sigma _1,\, ...,\, \sigma _n \in (0, \infty ) $
for some $n \in \mathbb{N} $ are denoted by $\mbox{$\mathfrak C$} ( \sigma _1,\,...,\, \sigma _n).$
We state some inequalities involving $s(x)$ or $x- \tau \, t\, e_1.$
\begin{lemma}[\mbox{\cite[Lemma 4.8]{DeKr2}}] \label{lemma2.10}
$s(x-y) ^{-1} \le \mbox{$\mathfrak C$} \, (1+|y|) \, s(x) ^{-1} $ for $x,y \in \mathbb{R}^3 $.
\end{lemma}
\begin{lemma}[\mbox{\cite[Lemma 2]{D1}}] \label{lemma2.20}
For $x \in \mathbb{R}^3 ,\; t \in (0, \infty ) , $ we have
\begin{eqnarray*} &&
|x- \tau \hspace{1pt}  t \hspace{1pt}  e_1|^2+t \ge \mbox{$\mathfrak C$} \hspace{1pt}
\bigl[\, \chi_{[0,1]}(|x|)\, (|x|^2+t) +
\chi_{(1, \infty )}(|x|)\, \bigl(\, |x|\, s(x)+t \,\bigr) \,\bigr]
.
\end{eqnarray*}
\end{lemma}
\begin{theorem}[\mbox{\cite[Theorem 2.19]{DKN2} }] \label{theorem2.50}
Let $R_1,\, R_2 \in (0, \infty ) $ with $R_1 < R_2,\; \nu \in (1, \infty ).$ Then
for
$y \in B_{R_2}^c,\; z \in B_{R_1}$,
\begin{eqnarray*}
\int_{ 0} ^{ \infty } (|y- \tau \, t\, e_1-e^{-t\, \Omega } \cdot z|^2 +t) ^{-\nu}\, dt
\le
\mbox{$\mathfrak C$} (R_1,R_2,\nu ) \, \bigl(\, |y|\, s(y) \,\bigr) ^{-\nu+1/2}.
\end{eqnarray*}
\end{theorem}
\begin{theorem} \label{theorem2.60}
Let $ R \in (0, \infty )$. Then for $ k \in \{0,\, 1\},\; x,y \in B_R$ with $x\ne y$,
\begin{eqnarray*}
\int_{ 0} ^{ \infty } \bigl(\, |x- \tau \, t \, e_1 - e^{-t \cdot \Omega }\cdot y|^2+t) ^{-3/2-k/2}\, dt
\le
C(R) \, |x-y|^{-1-k}.
\end{eqnarray*}
\end{theorem}
{\bf Proof:}
See the last part of the proof of \cite[Theorem 3.1]{DKN1}. Note that in \cite[(3.7)]{DKN1} it
should read $y+t\, U- e^{-t\, \Omega } \cdot z$ instead of $x$.
\hfill $\Box $
\vspace{1ex}

The next lemma is well known. It was already used in \cite{FH2}, for example.
For the convenience of the reader, we give a proof.
\begin{lemma} \label{lemma2.30}
Let $t \in \mathbb{R} .$ Then
$
e^{t \hspace{1pt}  \Omega }
=
\left( \begin{array}{rrr}
1 & 0  & 0  \\
0& \cos(t \hspace{1pt}  \varrho ) & - \sin(t \hspace{1pt}  \varrho )\\
0 & \sin (t \hspace{1pt}  \varrho ) & \cos( t \hspace{1pt}  \varrho )
\end{array}
\right).
$
\end{lemma}
{\bf Proof:}
Put
$
\widetilde{ I}
:=
\left( \begin{array}{rrr}
0 & 0  & 0  \\
0& 1 & 0\\
0 & 0 & 1
\end{array}
\right),
\;
I:=
\left( \begin{array}{rrr}
1 & 0  & 0  \\
0& 1 & 0\\
0 & 0 & 1
\end{array}
\right),
\;
A
:=
\left( \begin{array}{rrr}
0 & 0  & 0  \\
0& 0 & -1\\
0 & 1 & 0
\end{array}
\right).
$
It is easy to check that $ \Omega = \varrho \hspace{1pt}  A,\; A^2= -\widetilde{ I},\; A^3=-A$ and
$A^4 = \widetilde{ I}.$ Therefore
\begin{eqnarray*} &&
e^{t \hspace{1pt}  \Omega }
=
\displaystyle{\sum_{i=0} ^{ \infty }} ( \varrho \hspace{1pt}  t \hspace{1pt}  A)^i/i!
=\;
I \, +\,  \Bigl( \displaystyle{\sum_{i=1} ^{ \infty }}
( \varrho \hspace{1pt}  t)^{4 \hspace{1pt}  i} /(4 \hspace{1pt}  i)! \Bigr)   \hspace{1pt}  \widetilde{ I}
+
\Bigl( \displaystyle{\sum_{i=0} ^{ \infty }}
( \varrho \hspace{1pt}  t)^{4 \hspace{1pt}  i+1} / (4 \hspace{1pt}  i+1)! \Bigr)   \hspace{1pt}  A
\\ && \hspace{1em}
+
\Bigl( \displaystyle{\sum_{i=0} ^{ \infty }}
( \varrho \hspace{1pt}  t)^{4 \hspace{1pt}  i + 2} / (4 \hspace{1pt}  i +2)! \Bigr)
\hspace{1pt}  (-\widetilde{ I})
+
\Bigl( \displaystyle{\sum_{i=0} ^{ \infty }}
( \varrho \hspace{1pt}  t)^{4 \hspace{1pt}  i + 3} / (4 \hspace{1pt}  i+3)! \Bigr)   \hspace{1pt}  (-A)
\\ &&
=
I
 +
\Bigl(\displaystyle{ \sum_{j=1} ^{ \infty }} (-1)^j \hspace{1pt}
( \varrho \hspace{1pt}  t)^{2 \hspace{1pt}  j} / (2 \hspace{1pt}  j)! \Bigr)
\hspace{1pt}  \widetilde{ I}
 +
\Bigl(\displaystyle{ \sum_{j=0} ^{ \infty} } (-1)^j \hspace{1pt}
( \varrho \hspace{1pt}  t)^{2 \hspace{1pt}  j+1} / (2 \hspace{1pt}  j+1)! \Bigr) \hspace{1pt}  A
\\&&
=
I  +  (\cos( \varrho \hspace{1pt}  t) -1) \hspace{1pt}  \widetilde{ I} \; +\;
\sin( \varrho \hspace{1pt}  t) \hspace{1pt}  A.
\end{eqnarray*}
This implies the lemma.
\hfill $\Box $
\vspace{1ex}

Next we introduce some fundamental solutions. Put
\begin{eqnarray*}
N(x):=(4\, \pi\, |x|) ^{-1}
\quad \mbox{for}\;\;
x \in \mathbb{R}^3 \backslash\{0\}
\end{eqnarray*}
(''Newton potential'', fundamental solution of the Poisson equation in $\mathbb{R}^3 $),
\begin{eqnarray*}
\mbox{$\mathfrak O$}(x):=(4\, \pi\, |x|) ^{-1} \, e^{-\tau \,(|x|-x_1) /2}
\quad \mbox{for}\;\;
x \in \mathbb{R}^3 \backslash\{0\}
\end{eqnarray*}
(fundamental solution of the scalar Oseen equation $-\Delta v+ \tau \, \partial _1v=g$ in $\mathbb{R}^3 $),
\begin{eqnarray*}
\mbox{$\mathfrak O$}^{(\lambda )}(x):=(4\, \pi\, |x|) ^{-1} \, e^{-\sqrt{\lambda +\tau^2/4} \,|x|+\tau \,x_1 /2}
\quad \mbox{for}\;\;
x \in \mathbb{R}^3 \backslash\{0\} ,\; \lambda \in (0, \infty )
\end{eqnarray*}
(fundamental solution of the scalar Oseen resolvent equation $-\Delta v+ \tau \, \partial _1v + \lambda
\, v=g$ in $\mathbb{R}^3 $),
\begin{eqnarray*}
K(x,t):=(4\,\pi\, t) ^{-3/2} \, e^{-|x|^2/(4\, t)}
\quad \mbox{for}\;\;
x \in \mathbb{R}^3 ,\; t \in (0, \infty )
\end{eqnarray*}
(fundamental solution of the heat equation in $\mathbb{R}^3 $),
\begin{eqnarray*}&&
\psi(r):= \int_{ 0}^r (1-e^{-t})\, t ^{-1} \, dt\;\; (r \in \mathbb{R} ) ,
\quad
\Phi(x):=(4\, \pi\, \tau ) ^{-1} \, \psi \bigl(\, \tau \, (|x|-x_1)/2 \,\bigr)
\;\; (x \in \mathbb{R}^3 ),
\\&&
E _{jk} (x):=( \delta _{jk} \, \Delta - \partial _j \partial _k) \Phi(x),
\;\;
E_{4k}(x):=x_k\, (4 \, \pi \, |x| ^{3}) ^{-1}
\;\;
(x \in \mathbb{R}^3 \backslash\{0\} ,\; 1\le j,k\le 3)
\end{eqnarray*}
(fundamental solution of the Oseen system
(\ref{50}), with $(E _{jk} ) _{1\le j,k\le 3} $ the velocity part and $(E_{4k})_{1\le k\le 3}$ the pressure part).
We further define
\begin{eqnarray*}
F ^{( \lambda )}(\xi ):=(2\, \pi) ^{-3/2} \, ( \lambda + |\xi|^2 + i\, \tau \, \xi_1) ^{-1}
\quad \mbox{for}\;\;
\xi \in \mathbb{R}^3,\; \lambda \in (0, \infty )
\end{eqnarray*}
(Fourier transform of $\mbox{$\mathfrak O$}^{( \lambda )}$; see Theorem \ref{theorem3.10}).

We recall some basic properties of these functions, beginning with a classical result.
\begin{lemma} \label{lemma2.50}
Let $f \in \mbox{$\mathfrak S$} ( \mathbb{R}^3 )$ and put $F(x):=\int_{ \mathbb{R}^3 }N(x-y)\, f(y)\, dy$
for $x \in \mathbb{R}^3 .$ Then $F \in C ^{ \infty } ( \mathbb{R}^3 ) $
and $\partial ^{\alpha }F(x)= \int_{ \mathbb{R}^3 }N(x-y)\, \partial ^{\alpha }f(y)\, dy$
for $x \in \mathbb{R}^3 ,\; \alpha \in \mathbb{N} _0^3.$
\end{lemma}
\begin{lemma}[\mbox{\cite{Eidelman}}] \label{lemma2.60}
$K \in C ^{ \infty } \bigl(\, \mathbb{R}^3 \times (0, \infty ) \,\bigr) $ and
\begin{eqnarray*}
| \partial _t^l \partial _x^{\alpha }K(x,t)|
\le
\mbox{$\mathfrak C$} ( \alpha ,l) \, (|x|^2+t) ^{-3/2-| \alpha |/2-l}\, e^{-|x|^2/(8\, t)}
\end{eqnarray*}
for
$
x \in \mathbb{R}^3 ,\; t \in (0, \infty ) ,\; \alpha \in \mathbb{N} _0^3,\; l \in \mathbb{N} _0.
$
In particular $K(\; \cdot \; ,t) \in L^1( \mathbb{R}^3 ) \cap \mbox{$\mathfrak S$}( \mathbb{R}^3 )$ for
$t>0$.
\end{lemma}
\begin{theorem}[\mbox{\cite{KNP}}] \label{theorem2.70}
$E _{jk} \in C ^{ \infty } ( \mathbb{R}^3 \backslash\{0\} )$ and
\begin{eqnarray*}
| \partial ^{\alpha }E _{jk} (x)|\le \mbox{$\mathfrak C$} \, \bigl(\, |x|\, s(x) \,\bigr) ^{-1-| \alpha |/2}
\, \max\{1,\,|x|^{-| \alpha |/2}\}
\end{eqnarray*}
for
$
x \in \mathbb{R}^3 \backslash\{0\} ,\; 1\le j,k\le 3,\; \alpha \in \mathbb{N} _0^3\; \mbox{with}\;
| \alpha |\le 1.
$
\end{theorem}
As a consequence of Theorem \ref{theorem2.70}, we have $E _{jk} \in L^1_{loc}( \mathbb{R}^3 \backslash\{0\} )$
and $E _{jk} |B_1^c $ bounded $(1\le j,k\le 3)$. The same properties are obvious for $N,\,
\mbox{$\mathfrak O$}$ and $\mbox{$\mathfrak O$}^{( \lambda )}$.
Moreover $|\Phi (x)|\le \mbox{$\mathfrak C$} \, (1+|x|)\; (x \in \mathbb{R}^3 )$.
In view of these observations, the Fourier transforms of these functions will be considered as
tempered distributions (which, of course, will turn out to be represented by functions).
Following Solonnikov \cite[(40)]{Sol2}, we use Lemma \ref{lemma2.50} and \ref{lemma2.60} to
introduce the velocity part $(T _{jk} )_{1\le j,k\le 3}$ of a fundamental solution of the time-dependent
Stokes system, setting
$$
T _{jk} (x,t):= \delta _{jk} \, K(x,t) +\partial _j \partial _k \Bigl( \int_{ \mathbb{R}^3 }N(x-y)\, K(y,t)\, dy \Bigr)
\;\;
(x \in \mathbb{R}^3 ,\; t >0 ,\; 1\le j,k\le 3).
$$
\begin{lemma}[\mbox{\cite[Lemma 13]{Sol2}, \cite{Sol1}}] \label{lemma2.80}
$T _{jk} \in C ^{ \infty } \bigl(\, \mathbb{R}^3 \times (0, \infty ) \,\bigr) $ and
$$
| \partial _t^l \partial _x^{\alpha }T _{jk} (x,t)|\le \mbox{$\mathfrak C$} ( \alpha ,l)\,(|x|^2+t) ^{-3/2-| \alpha |/2-l}
$$
for
$
x \in \mathbb{R}^3 ,\; t \in (0, \infty ) ,\; 1\le j,k\le 3,\; \alpha \in \mathbb{N} _0^3,\;
l \in \mathbb{N} _0.
$
\end{lemma}
Lemma \ref{lemma2.80} yields that $T _{jk}(\; \cdot \; ,t) \in L^2 ( \mathbb{R}^3 ),$ but does not
imply $T _{jk}(\; \cdot \; ,t) \in L^1 ( \mathbb{R}^3 )\;\; (t>0)$. So the Fourier transform of this
function should be understood either as a transform of an $L^2$-function or as a tempered distribution.
For us it will be convenient to use the second possibility.
Put
\begin{eqnarray} \label{2.10}
\Gamma _{jk} (x,y,t):= T(x-\tau \, t\, e_1-e^{-t\, \Omega }\cdot y,\, t) \cdot e^{-t\, \Omega }
\quad \mbox{for}\;\;
x,y \in \mathbb{R}^3 ,\; t>0.
\end{eqnarray}
The matrix-valued function $( \Gamma _{jk} )_{1\le j,k\le 3}$ is the velocity part of a fundamental solution
to the time-dependent variant of the linearization (\ref{15}) of (\ref{10}). This fundamental solution was
constructed by Guenther, Thomann \cite{GT} via a procedure involving Kummer functions, an approach also
used in \cite{DKN1} -- \cite{DKN7}. However, Guenther, Thomann \cite[(3.9)]{GT} showed that
$\Gamma $ is given by (\ref{2.10}) as well, thus providing an access to this function which is more
convenient in many respects. For example, from Lemma \ref{lemma2.80} and (\ref{2.10}) we immediately
obtain
\begin{corollary} \label{corollary2.100}
Let $j,k \in \{1,\, 2,\, 3\} $. Then $\Gamma _{jk} \in C ^{ \infty } \bigl(\, \mathbb{R}^3
\times \mathbb{R}^3  \times (0, \infty )
\,\bigr) $ and
\begin{eqnarray*}
| \partial _x^{\alpha }\Gamma  _{jk} (x,y,t)|\le \mbox{$\mathfrak C$} ( \alpha )\,
(|x-\tau \, t\, e_1-e^{-t\, \Omega }\cdot y|^2+t) ^{-3/2-| \alpha |/2}
\end{eqnarray*}
for
$
x,y \in \mathbb{R}^3 ,\; t \in (0, \infty ) ,\;\alpha \in \mathbb{N} _0^3.
$
\end{corollary}
By Theorem \ref{theorem2.60} and Corollary \ref{corollary2.100}, we have
$
\int_{ 0} ^{ \infty } | \Gamma (x,y,t)|\, dt < \infty
$
for
$x,y \in \mathbb{R}^3$ with $x\ne y$, so we may define
\begin{eqnarray*}
Z _{jk} (x,y):= \int_{ 0} ^{ \infty }  \Gamma (x,y,t)\, dt
\quad \mbox{for}\;\;
x,y \in \mathbb{R}^3\; \mbox{with}\; x\ne y,\; 1\le j,k\le 3.
\end{eqnarray*}
This function $Z$ was introduced on \cite[p. 96]{GT} as the velocity part of a
fundamental solution to (\ref{15}). We collect the properties of $Z$ that will be needed
in what follows.
\begin{lemma}[\mbox{\cite[Lemma 2.15]{DKN2}}] \label{lemma2.100}
$
Z \in C^1 \bigl(\, ( \mathbb{R}^3 \times \mathbb{R}^3 ) \backslash \mbox{diag}\,( \mathbb{R}^3 \times \mathbb{R}^3 )
\,\bigr)^{3 \times 3} ,
\;\;
\partial x_lZ(x,y) = \int_{ 0} ^{ \infty } \partial x_l \Gamma (x,y,t)\, dt
$
for
$
x,y \in \mathbb{R}^3
$
with
$
x\ne y,\; 1\le l\le 3.
$
\end{lemma}
Note that due to Theorem \ref{theorem2.60} and Corollary \ref{corollary2.100}, we have
$
\int_{ 0} ^{ \infty } |\partial x_l \Gamma (x,y,t)|\, dt < \infty
$
for $x,\, y,\,l$ as in Lemma \ref{lemma2.100}.
\begin{corollary} \label{corollary2.110}
Let $R_1,\, R_2 \in (0, \infty ) $ with $R_1<R_2.$ Then
$$
| \partial _x^{\alpha }Z(x,y)|\le \mbox{$\mathfrak C$} (R_1,R_2) \, \bigl(\, |x|\, s(x) \,\bigr) ^{-1- |\alpha |/2}
\;\; \mbox{for}\;\;
x \in B_{R_2}^c,\; y \in B_{R_1},\; \alpha \in \mathbb{N} _0^3\; \mbox{with}\; | \alpha |\le 1.
$$
\end{corollary}
{\bf Proof:}
Lemma \ref{lemma2.100}, Corollary \ref{corollary2.100}, Theorem \ref{theorem2.50}. \hfill $\Box $
\begin{corollary} \label{corollary2.120}\label{corollary2.130}
The function $Z(\; \cdot \; ,0) $ belongs to $ C^1( \mathbb{R}^3 \backslash\{0\} )^{3 \times 3}.$

Let $S \in (0, \infty ) $. Then
$
| \partial _x^{\alpha }Z(x,0)|\le \mbox{$\mathfrak C$} (S) \, \bigl(\, |x|\, s(x) \,\bigr) ^{-1-| \alpha |/2}
$
for
$
x \in B_{S}^c,\; \alpha \in \mathbb{N} _0^3$ with $| \alpha |\le 1$.

Moreover
$|Z(x,0)|\le \mbox{$\mathfrak C$} \, |x|^{-1}$ for $x \in B_1 \backslash\{0\} $.
\end{corollary}
{\bf Proof:} The first part of the corollary follows from Lemma \ref{lemma2.100} and Corollary \ref{corollary2.110}.
The last estimate is a consequence of Corollary \ref{corollary2.100} and Theorem \ref{theorem2.60}.
\hfill $\Box $

\vspace{1ex}
Corollary \ref{corollary2.120} justifies to introduce the Fourier transform of $Z(\; \cdot \; ,0)$
in the sense of a tempered distribution.

\section{Statement of our main result.}
\setcounter{equation}{0}

It will be convenient to first recall the main result from \cite{DKN7}.
\begin{theorem}[\mbox{\cite[Theorem 3.1]{DKN7}}] \label{theorem2.10}
Let $\mbox{$\mathfrak D$} \subset \mathbb{R}^3 $ be open, $p \in (1, \infty ),\; f \in L^p( \mathbb{R}^3 )^3
$ with $supp(f)$ compact. Let $S_1 \in (0, \infty ) $ with $\overline{ \mbox{$\mathfrak D$} }\cup supp(f) \subset
B_{S_1}.$

Let $u \in L^6( \overline{ \mbox{$\mathfrak D$} }^c)^3\cap
W^{1,1}_{loc}( \overline{ \mbox{$\mathfrak D$} }^c)^3,\;
\pi \in L^2_{loc}( \overline{ \mbox{$\mathfrak D$} }^c)$ with
$\nabla u  \in L^2( \overline{ \mbox{$\mathfrak D$} }^c)^9,\; \mbox{div}\, u = 0$ and
\begin{eqnarray*} \label{ZU5}&&\hspace{-3em}
\int_{ \overline{ \mbox{$\mathfrak D$} }^c} \Bigl[\,
\nabla  u \cdot \nabla \varphi  \; +\;  \bigl(\, \tau \, \partial _1 u
 + \tau \, (u \cdot \nabla )u
- (\omega \times z) \cdot \nabla  u
+
\omega \times u \,\bigr) \cdot \varphi  -  \pi \, \mbox{div } \varphi \Bigr] \, dz
\\&&\nonumber \hspace{-3em}
= \int_{ \overline{ \mbox{$\mathfrak D$} }^c} f \cdot \varphi \, dz
\quad \mbox{for}\;\;
\varphi \in C ^{ \infty } _0( \overline{ \mbox{$\mathfrak D$} }^c)^3.
\end{eqnarray*}
(This means the pair $(u,\pi )$ is a Leray solution to (\ref{10}), (\ref{20}).)
Suppose in addition that
\begin{eqnarray} \label{ZU10}
\mbox{$\mathfrak D$} \; \mbox{is $C^2$-bounded},
\quad
u| \partial \mbox{$\mathfrak D$} \in W^{2-1/p,\, p}( \partial \mbox{$\mathfrak D$} )^3,
\quad
\pi | B_{S_1}\backslash \overline{ \mbox{$\mathfrak D$} } \in L^p(B_{S_1}\backslash \overline{ \mbox{$\mathfrak D$} }).
\end{eqnarray}
Let $n$ denote the outward unit normal to \mbox{$\mathfrak D$}, and define
\begin{eqnarray*} &&
\beta _k
:=
\int_{ \overline{ \mbox{$\mathfrak D$} }^c} f_k(y)\, dy
\\&& \hspace{1em}
+
\int_{ \partial \mbox{$\mathfrak D$} } \sum_{l = 1}^3 \bigl(\,
- \partial _lu_k(y) + \delta _{kl} \, \pi(y)
+
( \tau \, e_1- \omega \times y)_l \, u_k(y) - \tau \, (u_l\, u_k)(y)\,\bigr)
\, n_l(y)
\;do_y
\end{eqnarray*}
for $1\le k\le 3,$
\begin{eqnarray*} &&\hspace{-2em}
\mbox{$\mathfrak F$} _j(x)
:=
\int_{ \overline{ \mbox{$\mathfrak D$} }^c}
\Bigl[ \sum_{k = 1}^3   \bigl(\,  Z_{jk}
(x,y) - Z _{jk} (x,0) \,\bigr)
\, f_k(y)
-
\tau \cdot \sum_{k,l = 1}^3 Z_{jk}(x,y)\, (u_{l}\,\partial_l u_k) (y)
\Bigr]\, dy
\\&& \hspace{-2em}
+
\int_{ \partial \mbox{$\mathfrak D$} }
\sum_{k = 1}^3 \Bigl[
\bigl(\,  Z_{jk}(x,y) - Z_{jk} (x,0) \,\bigr)
\sum_{l = 1}^3 \bigl(\, - \partial _lu_k(y) + \delta _{kl} \, \pi(y)
+
( \tau \, e_1 - \omega \times y)_l \, u_k(y) \,\bigr)
\,
n_l(y)
\\&&
+ \bigl(\, E_{4j}(x- y)-E_{4j}(x) \,\bigr)
\, u_k(y) \, n_k(y)
\\&&
+
\sum_{l = 1}^3 \bigl(\,
\partial y_l Z_{jk}(x,y) \, (u_k \, n_l)(y)
+
\tau
Z_{jk}(x,0)\,
(u_l\, u_k\, n_l)(y) \,\bigr)
\Bigr] \, do_y
\end{eqnarray*}
for $x \in \overline{ B_{S_1}}^c,\; 1\le j\le 3.$
The preceding integrals are absolutely convergent.
Moreover $\mbox{$\mathfrak F$} \in C^1( \overline{ B_{S_1}}^c)^3$ and equation (\ref{70}) holds.
In addition, for any $S \in (S_1, \infty )$, there is a constant $C>0$ which depends on
$\tau ,\, \varrho ,\, S_1,\, S,\, f,\,u$ and $\pi$,
and which is such that
\begin{eqnarray*}
| \partial ^{\alpha }\mbox{$\mathfrak F$} (x)|
\le
C\, \bigl(\, |x|\, s(x) \,\bigr) ^{-3/2-| \alpha |/2} \, \ln (2+|x|)
\quad \mbox{for}\;\;
x \in \overline{ B_S}^c,\; \alpha \in \mathbb{N} _0^3\; \mbox{with}\; | \alpha |\le 1.
\end{eqnarray*}
\end{theorem}
In the preceding theorem, the coefficients $\beta _1,\, \beta _2,\, \beta _3$ and the function
\mbox{$\mathfrak F$} are defined in terms of integrals on $\partial \mbox{$\mathfrak D$} $ and
$\overline{ \mbox{$\mathfrak D$} }^c$. The integral over $\partial \mbox{$\mathfrak D$} $ may allow
to exploit boundary conditions verified
by $u$ or $\pi$. However, this way of introducing $\beta _1,\, \beta _2,\, \beta _3$ and \mbox{$\mathfrak F$}
requires the additional assumptions imposed on $\mbox{$\mathfrak D$} ,\, u$ and $\pi$ in (\ref{ZU10}).
If boundary conditions on $\partial \mbox{$\mathfrak D$} $ do not matter, we may drop (\ref{ZU10})
and consider $(u| \overline{ B_{S_0}}^c,\,\pi | \overline{ B_{S_0}}^c)$ instead of $(u,\pi )$,
where $S_0$ may be any number from $(0,S_1)$ with $\overline{ \mbox{$\mathfrak D$} }\cup supp(f) \subset
B_{S_0}.$
In view of interior regularity of $u$ and $\pi$, we may then define the coefficients $\beta _k$ and
the functions \mbox{$\mathfrak F$}  in terms of integrals over $\partial B_{S_0}$ and $\overline{ B_{S_0}}^c,$
obtaining an analogous result as the one in Theorem \ref{theorem2.10}, but with $B_{S_0}$ in the role
of \mbox{$\mathfrak D$}. Below we will present a variant of this idea which takes account of the additional
results in the work at hand (Corollary \ref{corollary2.10}).

The principal aim of this article consists in improving Theorem \ref{theorem2.10} in the way specified
in
\begin{theorem} \label{theorem2.20}
Let $\mbox{$\mathfrak D$} ,\, p,\, f,\, S_1,\, u,\,\pi $ satisfy the assumptions of Theorem \ref{theorem2.10},
including (\ref{ZU10}). Let $\beta _1,\, \beta _2,\, \beta _3$ and \mbox{$\mathfrak F$} be defined as
in Theorem \ref{theorem2.10}. Define the function \mbox{$\mathfrak G$} as in (\ref{ZU31}).

Then $\mbox{$\mathfrak G$} \in C^1( \overline{ B_{S_1}}^c)^3$, equation (\ref{80}) holds, and
for any $S \in (S_1, \infty ),$ there is a constant $C>0$ which depends on $\tau ,\, \varrho ,\, S_1,\, S,\,
f,\, u$ and $\pi$,
and which is such that
\begin{eqnarray*}
| \partial ^{\alpha }\mbox{$\mathfrak G$} (x)|
\le
C\, \bigl(\, |x|\, s(x) \,\bigr) ^{-3/2-| \alpha |/2} \, \ln (2+|x|)
\quad \mbox{for}\;\;
x \in \overline{ B_S}^c,\; \alpha \in \mathbb{N} _0^3\; \mbox{with}\; | \alpha |\le 1.
\end{eqnarray*}
\end{theorem}
We recall that the asymptotic behaviour of the function $E$
appearing in the leading term in (\ref{80}) is described in Theorem \ref{theorem2.70}.
As explained above, we may drop the assumptions in (\ref{ZU10}) if we replace $(u,\pi )$ by
 $(u| \overline{ B_{S_0}}^c,\,\pi | \overline{ B_{S_0}}^c)$,
with some suitably chosen number $S_0$. Here are the details.
\begin{corollary} \label{corollary2.10}
Take $\mbox{$\mathfrak D$} ,\, p,\, f,\, S_1,\, u,\,\pi $ as in Theorem \ref{theorem2.10},
but without requiring (\ref{ZU10}).
(This means that $(u,\, \pi )$ is only assumed to be a Leray solution of (\ref{10}), (\ref{20}).)
Put $\widetilde{ p}:=\min\{3/2,\, p\}.$

Then $u \in W^{2, \widetilde{ p}}_{loc}( \overline{  \mbox{$\mathfrak D$} }^c)^3$ and
$\pi \in W^{1, \widetilde{ p}}_{loc}( \overline{  \mbox{$\mathfrak D$} }^c)$.

Fix some number $S_0 \in (0,S_1) $ with
$\overline{ \mbox{$\mathfrak D$} }\cup supp(f) \subset  B_{S_0}$, and define
$\beta _1,\, \beta _2,\, \beta _3$ and \mbox{$\mathfrak F$} as in Theorem \ref{theorem2.10}, but
with \mbox{$\mathfrak D$} replaced by $B_{S_0}$, and $n(x)$ by $S_0 ^{-1} \, x$, for $x \in \partial B_{S_0}.$
Moreover, define \mbox{$\mathfrak G$} as in (\ref{ZU31}).

Then all the conclusions of Theorem \ref{theorem2.20} are valid.
\end{corollary}

\section{Some Fourier transforms.}
\setcounter{equation}{0}

In this section we show that $Z _{j1}(\; \cdot \; ,0)=E _{j1}$. To this end, we prove
that the Fourier transforms of these two functions coincide. To begin with,
we recall some well known facts about the Fourier transforms of some of the
fundamental solutions introduced in Section 2. Other intermediate results in this section
may also be well known (Corollary \ref{corollary3.20} for example), but since their proofs
are very short, we present them for completeness.
\begin{theorem} \label{theorem3.10}
For $\xi \in \mathbb{R}^3 \backslash\{0\} $, we have
$
\widehat{N}(\xi ) = (2\, \pi ) ^{-3/2} \, |\xi | ^{-2} .
$
If $f \in \mbox{$\mathfrak S$} ( \mathbb{R}^3 ),\; \xi \in \mathbb{R}^3 \backslash\{0\} $,
and if $F(x):=\int_{ \mathbb{R}^3 }N(x-y)\, f(y)\, dy$ for $x \in \mathbb{R}^3 $, then
$\widehat{F}(\xi )=
|\xi | ^{-2} \, \widehat{f}(\xi )
\; (\xi \in \mathbb{R}^3 \backslash\{0\} ).
$

Moreover
$
\bigl[\, K(\; \cdot \; ,t) \,\bigr]^{\wedge}(\xi ) = (2\, \pi ) ^{-3/2}\, e^{-t\,|\xi|^2}
$
for
$
\xi \in \mathbb{R}^3 ,\; t>0,
$
and
$
\widehat{ \mbox{$\mathfrak O$}^{( \lambda )}}(\xi ) = F^{( \lambda )}(\xi )
$
for
$
\xi \in \mathbb{R}^3 ,\; \lambda \in (0, \infty ) .
$
\end{theorem}
{\bf Proof:} For the first formula, the reader may consult \cite[Proposition 2.1.1]{Neri} and its proof.
The second equation follows from the first by a well known formula for the Fourier transform of a
convolution. As a direct reference we mention \cite[Lemma V.1.1]{Stein}. The third equation is well
known, and as concerns the forth, we refer to \cite[Theorem 2.1]{DV}.
\hfill $\Box $
\begin{corollary} \label{corollary3.10}
$\widehat{\mbox{$\mathfrak O$}} (\xi )= (2\,\pi) ^{-3/2} \, (i\, \tau \, \xi_1 +|\xi|^2) ^{-1} $
for $\xi \in \mathbb{R}^3 \backslash\{0\} $.
\end{corollary}
{\bf Proof:}
Let $\varphi \in \mbox{$\mathfrak S$} ( \mathbb{R}^3 ).$
For $n \in \mathbb{N} ,\; \xi \in \mathbb{R}^3 $, we have
$
| F ^{(1/n)}(\xi) \hspace{1pt}  \varphi (\xi)|
\le
\mbox{$\mathfrak C$} \,
|\xi| ^{-2} \hspace{1pt}  | \varphi (\xi )|.
$
But $ \int_{ \mathbb{R}^3 }|\xi| ^{-2}  \hspace{1pt}  | \varphi (\xi )|\; d \xi < \infty $,
because $\varphi $ is rapidly decreasing. Thus we get from Lebesgue's theorem
\begin{eqnarray*}
\mbox{$\mathfrak A$}
:= (2\, \pi ) ^{-3/2}\, \int_{ \mathbb{R}^3 }(i\, \tau \, \xi_1+|\xi|^2) ^{-1} \, \varphi (\xi )\, d\xi
=
\lim_{n\to \infty }\int_{ \mathbb{R}^3 } F ^{(1/n)}(\xi ) \hspace{1pt}    \varphi (\xi)\; d\xi.
\end{eqnarray*}
Due to the last equation in Theorem \ref{theorem3.10}, we may conclude
\begin{eqnarray} \label{C3.10.10}
\mbox{$\mathfrak A$}
=
\lim_{n\to \infty }
\int_{ \mathbb{R}^3 } \mbox{$\mathfrak O$}  ^{(1/n)}(x ) \hspace{1pt}  \widehat{\varphi} (x)\, dx.
\end{eqnarray}
But
$
|\mbox{$\mathfrak O$}  ^{(1/n)}(x ) \hspace{1pt}   \widehat{\varphi} (x)|
\le
\mbox{$\mathfrak C$} \hspace{1pt}  |x| ^{-1} \hspace{1pt}  | \widehat{\varphi }(x)|
$
for $n \in \mathbb{N} ,\; x \in \mathbb{R}^3 \backslash\{0\} ,$
with $ \int_{ \mathbb{R}^3 }|x| ^{-1} \hspace{1pt}  | \widehat{\varphi }(x)|\, dx < \infty $
because $\varphi $ hence $\widehat{ \varphi }$ is rapidly decreasing.
Thus equation (\ref{C3.10.10}) and Lebesgue's theorem yield
$
\mbox{$\mathfrak A$}
=
\int_{ \mathbb{R}^3 } \mbox{$\mathfrak O$}  (x ) \hspace{1pt}   \widehat{\varphi} (x)\, dx .
$
Since this is true for any $\varphi \in \mbox{$\mathfrak S$}  ( \mathbb{R}^3 ), $
the corollary follows.
\hfill $\Box $
\begin{corollary} \label{corollary3.20}
Let $t \in (0, \infty ) ,\; j,k \in \{1,\, 2,\, 3\} . $ Then
\begin{eqnarray*}
[T _{jk} (\;  \cdot \; , t)]^{\wedge}(\xi )
=
(2\,\pi ) ^{-3/2}  \, ( \delta _{jk} -\xi_j\,\xi_k \,|\xi| ^{-2} )\, e^{-t\, |\xi |^2}
\quad \mbox{for}\;\; \xi \in \mathbb{R}^3 \backslash\{0\} .
\end{eqnarray*}
\end{corollary}
{\bf Proof:}
We have $K(\; \cdot \; ,t) \in \mbox{$\mathfrak S$} ( \mathbb{R}^3 )$ (Lemma \ref{lemma2.60}). Therefore
by Lemma \ref{lemma2.50},
$
T _{jk} (x,t) = \delta _{jk} \, K(x,t)+\int_{ \mathbb{R}^3 }N(x-y)\, \partial _j \partial _kK(y,t)\, dy
\; (x \in \mathbb{R}^3 ).
$
Since $K(\; \cdot \; ,t) $ belongs to $ \mbox{$\mathfrak S$} ( \mathbb{R}^3 )$ hence
$\partial _j \partial _k K(\; \cdot \; ,t) $ does, too,
Corollary \ref{corollary3.20} follows from Theorem \ref{theorem3.10}.
\hfill $\Box $
\begin{corollary} \label{corollary3.30}
Let $j \in \{1,\, 2,\, 3\} ,\; t \in (0, \infty ) .$ Then
\begin{eqnarray*}
[\Gamma _{j1}(\; \cdot \; ,0,t)]^{\wedge}(\xi )
=
(2\,\pi ) ^{-3/2} \, ( \delta _{j1}-\xi_j\, \xi_1\,|\xi| ^{-2} )\, e^{-t\,(i\, \tau \, \xi_1+|\xi |^2)}
\quad \mbox{for}\;\;
\xi \in \mathbb{R}^3 \backslash\{0\} .
\end{eqnarray*}
\end{corollary}
{\bf Proof:}
By Lemma \ref{lemma2.30}, we have
$
\Gamma _{j1}(x,0,t) = \bigl(\, T(x-\tau \, t\, e_1,\, t)\, e^{-t\, \Omega } \,\bigr) _{j1}
=
T_{j1}(x-\tau \, t\, e_1,\, t),
$
so Corollary \ref{corollary3.30} follows from Corollary \ref{corollary3.20}.
\hfill $\Box $
\begin{corollary} \label{corollary3.40}
Let $j \in \{1,\, 2,\, 3\} ,\; t \in (0, \infty ) .$ Then
\begin{eqnarray*}
[Z _{j1}(\; \cdot \; ,0)]^{\wedge}(\xi )
=
(2\,\pi ) ^{-3/2} \, ( \delta _{j1}-\xi_j\, \xi_1\,|\xi| ^{-2} )\, (i\, \tau \, \xi_1+|\xi |^2) ^{-1}
\quad \mbox{for}\;\;
\xi \in \mathbb{R}^3 \backslash\{0\} .
\end{eqnarray*}
\end{corollary}
{\bf Proof:}
Let $\varphi \in \mbox{$\mathfrak S$} ( \mathbb{R}^3 ).$ With Corollary \ref{corollary2.100}, we get
\begin{eqnarray*}
A:= \int_{ \mathbb{R}^3 }\int_{ 0} ^{ \infty } | \Gamma _{j1}(x,0,t)\, \widehat{\varphi} (t)|\, dt\, dx
\le
\mbox{$\mathfrak C$} \, \int_{ \mathbb{R}^3 } \int_{ 0} ^{ \infty } (|x- \tau \, t\, e_1|^2+t) ^{-3/2}\, |
\widehat{\varphi }(x)|\, dt\, dx.
\end{eqnarray*}
Next we apply Lemma \ref{lemma2.20} to obtain
\begin{eqnarray*} && \hspace{-1em}
A
\le
\mbox{$\mathfrak C$} \, \Bigl( \int_{ \mathbb{R}^3 } \int_{ 1} ^{ \infty } t ^{-3/2}\, |\widehat{\varphi }(x)|\, dt\, dx
+
\int_{ B_1 } \int_{ 0} ^{1 } (|x|^2+t) ^{-3/2}\, |\widehat{\varphi }(x)|\, dt\, dx
\\&& \hspace{2em}
+
\int_{ B_1^c } \int_{ 0} ^{1 } (|x|\, s(x)+t) ^{-3/2}\, |\widehat{\varphi }(x)|\, dt\, dx \Bigr)
\\&&\hspace{-1em}
\le
\mbox{$\mathfrak C$} \, \Bigl( \int_{ \mathbb{R}^3 }  |\widehat{\varphi }(x)|\, dx
+
\int_{ B_1 } \int_{ 0} ^{1 } |x| ^{-3/2} \,t ^{-3/4}\, |\widehat{\varphi }(x)|\, dt\, dx
+
\int_{ B_1^c } \int_{ 0} ^{1 } (1+t) ^{-3/2}\, |\widehat{\varphi }(x)|\, dt\, dx \Bigr)
\\&&\hspace{-1em}
\le
\mbox{$\mathfrak C$} \, \int_{ \mathbb{R}^3 }  |\widehat{\varphi }(x)|\, (1+|x| ^{-3/2}) \, dx.
\end{eqnarray*}
Since $\varphi $ hence $\widehat{ \varphi }$ belongs to $\mbox{$\mathfrak S$} ( \mathbb{R}^3 ),$
we know that $\int_{ \mathbb{R}^3 }|\widehat{\varphi }(x)|\, (1+|x| ^{-3/2}) \, dx < \infty ,$
so $A < \infty $. Therefore we may apply Fubini's theorem, to obtain
\begin{eqnarray*} &&
\int_{ \mathbb{R}^3 }Z_{j1}(x,0)\, \widehat{\varphi }(x)\, dx
=
\int_{ 0} ^{ \infty } \int_{ \mathbb{R}^3 } \Gamma _{j1}(x,0,t)\, \widehat{\varphi }(x)\, dx\, dt,
\\&&
=
\int_{ 0} ^{ \infty }  \int_{ \mathbb{R}^3 }
(2\,\pi ) ^{-3/2} \, ( \delta _{j1}-\xi_j\, \xi_1\,|\xi| ^{-2} )\, e^{-t\,(i\, \tau \, \xi_1+|\xi |^2)}
\, \varphi (\xi )\, d\xi \, dt,
\end{eqnarray*}
where the last equation follows from Corollary \ref{corollary3.30}. But
\begin{eqnarray*} &&
\int_{ 0} ^{ \infty }  \int_{ \mathbb{R}^3 }
| ( \delta _{j1}-\xi_j\, \xi_1\,|\xi| ^{-2} )\, e^{-t\,(i\, \tau \, \xi_1+|\xi |^2)}
\, \varphi (\xi )|\, d\xi \, dt
\le
\mbox{$\mathfrak C$} \,
\int_{ 0} ^{ \infty }  \int_{ \mathbb{R}^3 }  e^{-t\,|\xi |^2} \,| \varphi (\xi )|\, d\xi \, dt
\\&&
\le
\mbox{$\mathfrak C$} \, \int_{ \mathbb{R}^3 }  |\xi | ^{-2}  \,| \varphi (\xi )|\, d\xi < \infty ,
\end{eqnarray*}
with the last relation holding due to the assumption $\varphi \in \mbox{$\mathfrak S$} ( \mathbb{R}^3 )$.
Thus
we may use Fubini's theorem, arriving at the equation
\begin{eqnarray*}
\int_{ \mathbb{R}^3 } Z _{j1}(x,0)\, \widehat{\varphi }(x)\, dx
=
\int_{ \mathbb{R}^3 }
(2\,\pi ) ^{-3/2} \, ( \delta _{j1}-\xi_j\, \xi_1\,|\xi| ^{-2} )\, (i\, \tau \, \xi_1+|\xi |^2) ^{-1}
\, \varphi (\xi )\, d\xi.
\end{eqnarray*}
This proves Corollary \ref{corollary3.40}.
\hfill $\Box $
\begin{theorem} \label{theorem3.20}
Let $j,k \in \{1,\, 2,\, 3\} .$ Then for $\xi \in \mathbb{R}^3 \backslash\{0\} $,
\begin{eqnarray*}
\widehat{E} _{jk} (\xi ) = (2\, \pi ) ^{-3/2}\,
( \delta _{j1}-\xi_j\, \xi_1\,|\xi| ^{-2} )\, (i\, \tau \, \xi_1+|\xi |^2) ^{-1}
.
\end{eqnarray*}
\end{theorem}
{\bf Proof:}
For $x \in \mathbb{R}^3 \backslash\{0\} ,$ we find
\begin{eqnarray*} && \hspace{-2em}
\partial _1\Phi(x)
=
(4\,\pi\, \tau ) ^{-1} \,\psi ^{\prime} \bigl(\, \tau \,(|x|-x_1)/2 \,\bigr) \, \tau \, (x_1/|x|-1)/2
= (4\,\pi\, \tau\, |x| ) ^{-1} \,(e^{-\tau \,(|x|-x_1)/2}-1)
\\&&\hspace{-2em}
=
\tau ^{-1} \bigl(\, \mbox{$\mathfrak O$}(x)-N(x) \,\bigr) .
\end{eqnarray*}
Hence with Corollary \ref{corollary3.10} and Theorem \ref{theorem3.10}, for $\xi \in \mathbb{R}^3 \backslash\{0\} $,
\begin{eqnarray*} &&
i\, \xi_1\, \widehat{\Phi}(\xi )
=
\widehat{\partial _1\Phi}( \xi )
=
\tau ^{-1} \, (2\,\pi ) ^{-3/2}\, \bigl(\, (i\, \tau \,\xi_1+|\xi|^2) ^{-1} -|\xi| ^{-2} \,\bigr)
\\&&
=
-i\, (2\,\pi ) ^{-3/2}\, \xi_1\, \bigl(\, (i\, \tau \,\xi_1+|\xi|^2)\, |\xi|^2  \,\bigr) ^{-1} .
\end{eqnarray*}
As a consequence
$
\widehat{\Phi}(\xi )
=
-(2\,\pi ) ^{-3/2}\, \xi_1\, \bigl(\, (i\, \tau \,\xi_1+|\xi|^2)\, |\xi|^2  \,\bigr) ^{-1} ,
$
so the theorem follows by the definition of $E _{jk} $.
\hfill $\Box $

\vspace{1ex}
Theorem \ref{theorem3.20} may be deduced also from the results in \cite[Chapter VII]{Ga1}. In fact,
it is shown in \cite[Section VII.3]{Ga1} that the convolution $\mbox{$\mathfrak O$}*f$, for
$f \in C ^{ \infty } _0( \mathbb{R}^3 )^3,$ belongs to $C ^{ \infty } ( \mathbb{R}^3 )^3 $ and is the
velocity part of a solution to the Oseen system (\ref{50}) in $\mathbb{R}^3 $. On the other hand, by
\cite[Section VII.4]{Ga1}, the inverse Fourier transform of the function
$
(2\, \pi ) ^{-3/2}\,
( \delta _{j1}-\xi_j\, \xi_1\,|\xi| ^{-2} )\, (i\, \tau \, \xi_1+|\xi |^2) ^{-1} \, \widehat{f}(\xi )$
also solves (\ref{50}) in $\mathbb{R}^3 $, and belongs to certain Sobolev spaces. A uniqueness result
yields that the two solutions coincide, implying Theorem \ref{theorem3.20}.
However, we prefer to carry out a direct proof of this theorem, instead of relying on the rather lengthy
theory in \cite[Chapter VII]{Ga1}, which in fact yields much stronger results, not needed here, than
Theorem \ref{theorem3.20}.

Combining Theorem \ref{theorem3.20} and Corollary \ref{corollary3.40}, we arrive at the main result
of this section.
\begin{corollary} \label{corollary3.50}
$Z_{j1}(\; \cdot \; ,0)=E_{j1}$ for $1\le j\le 3$.
\end{corollary}

\section{Proof of Theorem \ref{theorem2.20} and Corollary \ref{corollary2.10}.}
\setcounter{equation}{0}

We first show that in the case $k \in \{2,\, 3\} ,$ the function $\partial ^{\alpha }_{jk}Z (\; \cdot \; ,0)$
decays faster for $|x|\to \infty $ than indicated by Corollary \ref{corollary2.120}.
\begin{theorem} \label{theorem4.10}
Let $S \in [2\, \tau \,\pi /| \varrho |,\, \infty ).$ Then
$
| \partial ^{\alpha }_xZ _{jk} (x,0)|\le \mbox{$\mathfrak C$} (S) \, \bigl(\, |x|\, s(x) \,\bigr) ^{-3/2-| \alpha |/2}
$
for
$
x \in B^c_{S+\tau \,\pi /| \varrho|},
\;
\alpha \in \mathbb{N} _0^3
$
with
$| \alpha |\le 1,\; j \in \{1,\, 2,\, 3\} ,\; k \in \{2,\, 3\} .$
\end{theorem}
{\bf Proof:}
Take $x,\, \alpha ,\, j,\, k$ as in the theorem.
We get with Lemma \ref{lemma2.100} that
\begin{eqnarray*}
\partial ^{\alpha }_xZ _{jk} (x,0)
=
\int_{ 0} ^{ \infty } \partial _x^{\alpha }\Gamma _{jk} (x,0,t)\, dt
=
\int_{ 0} ^{ \infty } \bigl[\, \partial _x^{\alpha }T(x-\tau \, t\, e_1,\, t) \cdot e^{-t\, \Omega } \,\bigr] _{jk} \, dt,
\end{eqnarray*}
so with Lemma \ref{lemma2.30} in the case $k=2$,
\begin{eqnarray} \label{T4.10.10}
\partial ^{\alpha }_xZ _{jk} (x,0)
=
\int_{ 0} ^{ \infty } \bigl(\, \partial _x^{\alpha }T_{j2}(x-\tau \, t\, e_1,\, t)\, \cos( \varrho \, t)
-\partial _x^{\alpha }T_{j3}(x-\tau \, t\, e_1,\, t)\, \sin( \varrho \, t) \,\bigr) \, dt,
\end{eqnarray}
with a similar formula in the case $k=3$. Let $\sigma : \mathbb{R} \mapsto \mathbb{R} $ be defined by
either $\sigma (t):=\cos ( \varrho \, t)$ for $t \in \mathbb{R} ,$ or by
$\sigma (t):=\sin ( \varrho \, t)$ for $t \in \mathbb{R} .$ Let $m \in \{1,\, 2,\, 3\} $. Then
\begin{eqnarray*} &&\hspace{-1em}
\int_{ 0} ^{ \infty } \partial _x^{\alpha }T_{jm}(x-\tau \, t\, e_1,\, t)\, \sigma (t)\, dt
\\&&\hspace{-1em}
=
\sum_{n = 0}^{ \infty } \int_{ n\,\pi/|  \varrho |} ^{ (n+1)\,\pi/|  \varrho |}
\Bigl(  \partial _x^{\alpha }T_{jm}(x-\tau \, t\, e_1,\, t)
-
\partial _x^{\alpha }T_{jm}(x-\tau \, (t+\pi/| \varrho |)\, e_1,\, t+\pi /| \varrho |)
\Bigr)  \, \sigma (t)\, dt
\\&&\hspace{-1em}
=
\sum_{n = 0}^{ \infty } \int_{ n\,\pi/|  \varrho |} ^{ (n+1)\,\pi/|  \varrho |}\int_{ 0}^1
(-\tau \,\partial _x^{\alpha +e_1 } +\partial _x^{\alpha }\partial _4)
T_{jm}(x-\tau \, (t+ \vartheta \,\pi/| \varrho |)\, e_1,\, t+\vartheta \,\pi /| \varrho |)
\\&&\hspace{10em}
\cdot  (-\pi/| \varrho |) \, \sigma (t)\, d \vartheta \, dt.
\end{eqnarray*}
Therefore by Lemma \ref{lemma2.80},
\begin{eqnarray*} && \hspace{-1em}
A:=
\Bigl| \int_{ 0} ^{ \infty } \partial _x^{\alpha }T_{j,m}(x-\tau \, t\, e_1,\, t)\, \sigma (t)\, dt \Bigr|
\\&&\hspace{-1em}
\le
\mbox{$\mathfrak C$} \, \sum_{n = 0}^{ \infty } \sum_{m = 1}^2 \int_{ n\,\pi/|  \varrho |} ^{ (n+1)\,\pi/|  \varrho |}\int_{ 0}^1
\bigl(\, |x-\tau \, (t+ \vartheta \,\pi/| \varrho |)\, e_1|^2 + t+\vartheta \,\pi /| \varrho |
\,\bigr) ^{-3/2-| \alpha |/2-m/2}
\, d \vartheta \, dt
\\&&\hspace{-1em}
\le
\mbox{$\mathfrak C$} \, \sum_{m = 1}^2 \int_{ 0}^1
\int_{0} ^{ \infty }
\bigl(\, |x-(\tau \,  \vartheta \,\pi/| \varrho |)\, e_1 - \tau \, t\, e_1|^2 + t \,\bigr) ^{-3/2-| \alpha |/2-m/2}
\, dt\, d \vartheta .
\end{eqnarray*}
Since
$x \in B^c_{S+\tau \,\pi /| \varrho|},$
we have
$
|x-( \tau \, \vartheta \, \pi /| \varrho |)\,e_1|\ge S
$
for $\vartheta \in [0,1] $, so we may apply Theorem \ref{theorem2.50}
with $z=0,\; R_2=S,\; R_1=S/2,\; y=x-( \tau \, \vartheta \, \pi/| \varrho |)\, e_1,\;
\nu = 3/2+| \alpha |/2+l/2,$ to obtain
\begin{eqnarray} \label{T4.10.20}
A
\le
\mbox{$\mathfrak C$}(S) \, \sum_{m = 1}^2 \int_{ 0}^1
\bigl[\, |x-(\tau \,  \vartheta \,\pi/| \varrho |)\, e_1|\:
s \bigl(\, x-(\tau \,  \vartheta \,\pi/| \varrho |)\, e_1 \,\bigr)  \,\bigr]
^{-1-| \alpha |/2-m/2} \, d \vartheta .
\end{eqnarray}
But for $\vartheta \in [0,1],$ we have
$
|x-(\tau \,  \vartheta \,\pi/| \varrho |)\, e_1|\ge |x|/2 +S/2-\tau \, \vartheta \,\pi/| \varrho |
\ge |x|/2,
$
where the last inequality holds because $S\ge 2\, \tau \, \pi/| \varrho|.$
Moreover we get from from Lemma \ref{lemma2.10} that
$
s \bigl(\, x-(\tau \,  \vartheta \,\pi/| \varrho |)\, e_1 \,\bigr)  ^{-1} \le \mbox{$\mathfrak C$} \, s(x) ^{-1}
$
for
$
\vartheta \in [0,1].
$
Therefore from (\ref{T4.10.20}),
\begin{eqnarray*}
A\le
\mbox{$\mathfrak C$}(S) \, \sum_{m = 1}^2 \bigl(\, |x|\, s(x) \,\bigr) ^{-1-| \alpha |/2-m/2}
\le
\mbox{$\mathfrak C$}(S) \,  \bigl(\, |x|\, s(x) \,\bigr) ^{-3/2-| \alpha |/2}.
\end{eqnarray*}
Theorem \ref{theorem4.10} follows with equation (\ref{T4.10.10}) and its analogue for $k=3$.
\hfill $\Box $
\begin{corollary} \label{corollary5.10}
Let $S \in (0, \infty ) $. Then
$
| \partial ^{\alpha }_xZ _{jk} (x,0)|\le \mbox{$\mathfrak C$} (S) \, \bigl(\, |x|\, s(x) \,\bigr) ^{-3/2-| \alpha |/2}
$
for
$
x \in B^c_{S}
$
and for
$ \alpha ,\, j,\, k$ as in Theorem \ref{theorem4.10}.
\end{corollary}
{\bf Proof:}
Let $x \in B_S^c$, and take $ \alpha ,\, j,\, k$ as in Theorem \ref{theorem4.10}.
By Corollary \ref{corollary2.120}, we have
$
| \partial ^{\alpha }_xZ _{jk} (x,0)|\le \mbox{$\mathfrak C$} (S) \, \bigl(\, |x|\, s(x) \,\bigr) ^{-1-| \alpha |/2}
$.

Suppose that
$S\ge 2\,\tau \,\pi /| \varrho |$.
Then we distinguish the cases
$x \in B_{S+\tau \,\pi /| \varrho |}^c$
and
$x \in B_{S+\tau \,\pi /| \varrho |}\backslash B_S$.
If $x \in B_{S+\tau \,\pi /| \varrho |}^c$,
the looked-for inequality follows from Theorem \ref{theorem4.10}.
In the second case, we observe that
$1\le (S+\tau \,\pi /| \varrho |)\, |x| ^{-1} ,$
so the inequality claimed in Corollary \ref{corollary5.10} may be deduced from the
estimate stated at the beginning of this proof.

Now suppose that
$S < 2\,\tau \,\pi /| \varrho |$,
Then we use that either
$x \in B_{3\,\tau \,\pi /| \varrho |}^c$
or
$x \in B_{3\,\tau \,\pi /| \varrho |}\backslash B_S.$
If $x \in B_{3\,\tau \,\pi /| \varrho |}^c$,
the looked-for inequality follows from Theorem \ref{theorem4.10} with
$2\,\tau \,\pi /| \varrho |$ in the place of $S$.
In the case
$x \in B_{3\,\tau \,\pi /| \varrho |}\backslash B_S$,
we use the relation
$1\le (3\,\tau \,\pi /| \varrho |)\, |x| ^{-1} $
and again the estimate from the beginning of the proof, once more
obtaining an upper bound
$
\mbox{$\mathfrak C$}(S) \,  \bigl(\, |x|\, s(x) \,\bigr) ^{-3/2-| \alpha |/2}
$
for
$
| \partial ^{\alpha }_xZ _{jk} (x,0)|,
$
as stated in Corollary \ref{corollary5.10}.
\hfill $\Box $

\vspace{1ex}
The proofs of Theorem \ref{theorem2.20} and Corollary \ref{corollary2.10} are now obvious.

{\bf Proof of Theorem \ref{theorem2.20}:}
Combine Theorem \ref{theorem2.10}, Corollary \ref{corollary3.50} and \ref{corollary5.10}.
\hfill $\Box $
\vspace{1ex}

{\bf Proof of Corollary\ref{corollary2.10}}:
From interior regularity of solutions to the Stokes system (\cite[Theorem IV.4.1]{Ga1}) and
the assumption $f \in L^p( \mathbb{R}^3 )^3,$ we may conclude that
$u \in W^{2, \widetilde{ p}}_{loc}( \overline{  \mbox{$\mathfrak D$} }^c)^3$ and
$\pi \in W^{1, \widetilde{ p}}_{loc}( \overline{  \mbox{$\mathfrak D$} }^c)$.
More details about this conclusion may be found in the proof of \cite[Theorem 5.5]{DKN2}.
It follows that
$
u| \partial B_{S_0} \in W^{2-1/ \widetilde{ p},\, \widetilde{ p}}( \partial B_{S_0})^3
$
and
$
\pi | B_{R}\backslash \overline{B_{S_0} } \in L^{\widetilde{ p}}(B_{R}\backslash \overline{B_{S_0} })
$
for any
$
R \in (S_0, \infty ).
$
Now we may apply Theorem \ref{theorem2.20} with $\mbox{$\mathfrak D$} ,\, f,\,u,\, \pi$
replaced by $B_{S_0},\, f| \overline{ B_{S_0}}^c,\, u| \overline{ B_{S_0}}^c$ and
$\pi | \overline{ B_{S_0}}^c$, respectively. Corollary \ref{corollary2.10} then follows
from Theorem \ref{theorem2.20}.
\hfill $\Box $

\vspace{2ex}
{\bf Acknowledgements.}  {\it The work of \v S\' arka Ne\v casov\' a
acknowledges the support of the GA\v CR (Czech Science Foundation)
project P201-13-00522S in the framework of RVO: 67985840. \v S.N. would like to thank for fruitful discussions with R. Guenther.}

\end{document}